\documentclass[11pt,fullpage,reqno]{article}
\usepackage{mathrsfs}%                                          *
\usepackage{pifont}%                                            *
\usepackage{amsmath}%                                           *
\usepackage{amsthm}%                                            *
\usepackage{txfonts}%                                           *
 \usepackage{geometry}%                                         *
\usepackage{latexsym}%                                          *
\usepackage{amssymb}%                                           *
\usepackage{graphicx}%                                          *
\usepackage{geometry}%                                          * Please do not change any words here.
\usepackage{ifthen}%                                            *
\usepackage{fancyhdr}%                                          *
\usepackage{afterpage}%                                         *
\usepackage{xcolor}%   
\usepackage{ifthen}%                                            *
\newlength{\defbaselineskip}
\newcommand{\setlinespacing}[2]%
           {\setlength{\baselineskip}{#1 \defbaselineskip}}

\begin{document}
\begin{center}
{\Large {\bf Estimation of common change point and isolation of changed panels after sequential detection}}\\

{\it Yanhong Wu \footnote{Email: ywu1@csustan.edu }} \\

California State University Stanislaus, Turlock CA 95382, USA 
\end{center}

{\small 
\noindent{\bf ABSTRACT:}
Quick detection of common changes is critical in sequential monitoring of multi-stream data where a common change is referred as a change that only occurs in a portion of panels.  After a common change is detected by using a combined CUSUM-SR procedure, we first study the joint distribution for values of the CUSUM process and the estimated delay detection time for the unchanged panels.  The BH method by using the asymptotic exponential property for the CUSUM process is developed to isolate the changed panels with the control on FDR. The common change point is then estimated based on the isolated changed panels. Simulation results show  that the proposed method can also control the FNR by properly selecting FDR. \\

}

\section{Introduction}
Detecting a common change in multi-stream data or panel data is critical on sequential change-point detection problem. Here, a common change is referred as a change that may occur only in a portion of the N panels; usually caused by external sources. In contrast, the traditional change-point detection is focused on a single sequence (individual panels) where the change is typically caused from internal sources.
Several typical detection procedures have been discussed and extended; see Xie and Siegmund (2013), Mei (2013), and Tartakovsky and Veeravalli(2008).
Chan (2017) discussed the optimality of detection procedures.

Wu (2019) proposed a combined SR-CUSUM procedure that uses the sum of N Shiryayev-Roberts processes to detect the common change, while the N individual CUSUM processes are used to isolate the changed panels and estimate the change point. The alarming limit $B$ is chosen such that the average in-control run length is equal to a designated value. For convenience of discussion, we shall focus on the normal case.

Assume there are $N$ independent panels and in panel $i$, the observations $\{ X_j^{(i)}\}$ follow $N(0,1)$ for $j \leq \nu $ and $N(\mu, 1)$ ($\mu >0$) for $j > \nu$ if a change occurs at $\nu $ in this panel. Suppose a change may only occur in $K$ of the $N$ panels, called common change. The $N$ panels can be assumed following a mixture model with probability $p=K/N$ of change in each panel. We shall select the same reference parameter $\delta $ for $\mu$ for all panels.

For $R_0(i)=0$, define
\[
R_t(i) = (1+R_{t-1}(i)) e^{\delta X_j^{(i)} -\delta^2/2 } = \sum_{k=1}^t e^ {\sum_{k+1}^t (\delta X_j^{(i)} -\delta^2/2)}
\]
as the Shiryayev-Roberts process for the $i^{th}$ panel and $R_t =\sum_{i=1}^N R_t(i) $ for $R_0 =0$.

An alarm will be raised at the stopping time
\[
\tau= \inf\{t>0:  R_t > B\},
\]
where $B$ is chosen such that the $ARL_0 $ is equal to the designated value.

In the normal case, when $\delta $ is small, $B$ can be designed by using the following simple approximation (Pollak (1987)):
\[
ARL_0 =E_0 \tau \approx \frac{B}{N} e^{\rho \delta} ,
\]
where $\rho \approx 0.5826 $.

For example, for $\delta =0.5$, $N=100$, $B= 74729.5$ and $373645.7$ corresponding to $ARL_0 =1000$ and $5000$, respectively.
Further properties on the average run lengths are referred to Wu (2018b).
Comparison with other procedures as shown in Wu (2018b) demonstrated that the proposed procedure is very competitive when the proportion is small.

To isolate the changed panels and estimate the common change point, the combined SR-CUSUM procedure calculates the $N$ CUSUM processes recursively as
\[
T_{t}(i) =\max (0, T_{t-1}(i) +  X_t^{(i)} -\delta/2).
\]
and at the alarm time $\tau$ the change-point for the $i^{th}$ panel is estimated as the last zero point of $T_t (i)$,
\[
\hat{\nu}_i = \max \{ t<\tau: T_{t}(i) = 0 \},
\]
for $i=1,..., N$, which is indeed the MLE when $\mu =\delta$. When $\mu$ is unknown, it can be estimated as
\[
\hat{\mu}_i = \frac{T_{\tau} (i)}{\tau-\hat{\nu}_i }+ \frac{\delta}{2}.
\]

Apparently, to isolate the "true" changed panels, both the change-point estimation (or estimated delay detection time $\tau-\hat{\nu}_i $ after the change-point estimation) and estimated strength of signals (or $T_{\tau}(i)$) provide related information. Here we propose a BH-type procedure to control the FDR.
In Section 2, we shall first study the corresponding continuous time Brownian motion model and derive the exact null joint distribution for $T_{\tau}(i)$ and $\tau-\hat{\nu}_i$ for the unchanged panels. The marginal moments and covariance shows that they are highly correlated. When $\delta $ is small, we extend the results to the discrete time model in Section 3. In Section 4, we propose to use the approximate null distribution for $T_{\tau}(i)$ to form a BH-type procedure to isolate the changed panel by controlling FDR. The isolated changed panels are then used to estimate the common change point.
Simulation studies for the FDR, FNR, and biases of estimated common change-point in several typical cases show the proposed method works quite well. The results also help us to select the proper FDR in order to balance FNR.

\section{Null Distribution under continuous time model}
We assume that a common change is detected (B is large) and the change occurs far away from the beginning. For those unchanged panels, at the detection time, by looking backward at each CUSUM process and using its strong Markov property, we can see that for each $i$, $\tau-\hat{\nu}_i$ and $T_{\tau}(i)$ are approximately equivalent in distribution to the maximum point $\sigma_M$ and the maximum value $M$ for a normal random walk $S_n=
\sum_{i=1}^n (X_i -\delta /2)$ for $S_0=0$ with drift $-\delta /2$ and variance 1 where
\[
M =\sup_{0\leq n < \infty} S_n ~~~~~~and~~~~~~~\sigma_M =argsup_{0\leq n < \infty} S_n .
\]

Under the continuous time model, we shall denote $\{W_t\} $ for $t\geq 0$ as a Brownian motion and $P(.)$ as its corresponding probability measure with drift $-\delta/2$ and $P^*(.)$ as the probability measure when the drift is $\delta/2$.
Denote
\[
M=\sup_{0\leq t< \infty} W_t ~~~~~~~ and ~~~~~~\sigma_M =argsup_{0\leq t <\infty} W_t.
\]
For an independent copy $W_t'$ of $W_t$, we denote $M'=\sup_{0\leq t<\infty}  W_t'$. The following theorem gives the joint distribution of $(\sigma_M, M)$ and its proof is given in the appendix.

\noindent{\bf Theorem 1.} {\it
\[
P[\sigma_M<t, M>x]=
\]
\[
P^*[\sup_{0\leq s<t} W_s >M']-E[e^{\delta M'} \Phi (-\frac{x+M'}{\sqrt{t}} -\frac{\delta}{2}\sqrt{t} )]
-E[e^{\delta x} \Phi (-\frac{x+M'}{\sqrt{t}} +\frac{\delta}{2}\sqrt{t} )].
\]
}

 Note that by letting $t \rightarrow \infty$, we see
 \[
 P[M>x] =e^{-\delta x}.
 \]
By taking derivative with respect to x, we get
 \begin{eqnarray*}
 P[\sigma_M <t, M\in dx] &=& E[\frac{1}{\sqrt{t}} e^{M' \delta} \phi (\frac{x+M'}{\sqrt{t}} +\frac{\delta}{2}\sqrt{t} )
 +\delta e^{-\delta x} \Phi (-\frac{x+M'}{\sqrt{t}} +\frac{\delta}{2}\sqrt{t} )] \\
 & & +\frac{1}{\sqrt{t}} e^{-\delta x} \phi (-\frac{x+M'}{\sqrt{t}} +\frac{\delta}{2}\sqrt{t} ) \\
 &=& E[\delta e^{-\delta x} \Phi (-\frac{x+M'}{\sqrt{t}} +\frac{\delta}{2}\sqrt{t} ) +\frac{2}{\sqrt{t}} e^{-\delta x} \phi (-\frac{x+M'}{\sqrt{t}} +\frac{\delta}{2}\sqrt{t} )], \\
 \end{eqnarray*}
 since
 \[
 e^{M' \delta} \phi (\frac{x+M'}{\sqrt{t}} +\frac{\delta}{2}\sqrt{t} ) =e^{-\delta x} \phi (-\frac{x+M'}{\sqrt{t}} +\frac{\delta}{2}\sqrt{t} )].
 \]
 Thus,
 \[
 P[\sigma_M<t|M=x]= E[ \Phi (-\frac{x+M'}{\sqrt{t}} +\frac{\delta}{2}\sqrt{t} ) +\frac{2}{\delta \sqrt{t}} \phi (-\frac{x+M'}{\sqrt{t}}+\frac{\delta}{2}\sqrt{t} )].
 \]
 The following theorem shows that the conditional distribution of $\sigma_M$ given $M=x$ is actually inverse Gaussian $IG( \frac{x}{\delta/2}, x^2 )$ under $P^*(.)$ and its proof is given in the appendix.

 \noindent{\bf Theorem 2.} 
 {\it
 \begin{eqnarray*}
  P[\sigma_M<t|M=x] &=& P^*[\tau_x < t ] \\
   &=& e^{\delta x} \Phi(-\frac{x}{\sqrt{t}} -\frac{\delta}{2} \sqrt{t} ) +\Phi (-\frac{x}{\sqrt{t}} +\frac{\delta}{2} \sqrt{t} ),
  \end{eqnarray*}
where $\tau_x =\inf\{t \geq>0: W_t >x \}$ and the conditional density function of $\sigma_M $ given $M=x$ is
\[
f_{\sigma_M|M} (t|x) = \frac{x}{t^{3/2} } \phi (-\frac{x}{\sqrt{t}} +\frac{\delta}{2} \sqrt{t}),
\]
and the joint density function of $(\sigma_M, M)$ is
\[
f(t, x) =\frac{\delta x}{t^{3/2}} \phi (\frac{x}{\sqrt{t}} +\frac{\delta}{2} \sqrt{t}) .
\]
}

By integrating $f(t,x)$ with respect to $x$, we have:

\noindent{\bf Corollary 1.}  {\it The marginal density function and cdf of $\sigma_M$  are given by
\begin{eqnarray*}
f_{\sigma_M}(t) &=& \frac{\delta}{\sqrt{t}}  \int_0^{\infty} (1-\Phi(y+\frac{\delta}{2} \sqrt{t})) dy \\
&=& \frac{\delta}{\sqrt{t}} \phi (\frac{\delta}{2}\sqrt{t})-\frac{\delta^2}{2} (1-\Phi (\frac{\delta}{2} \sqrt{t} )).
\end{eqnarray*}
\[
P[\sigma_M>t] = (4+\frac{\delta^2}{2}t)(1-\Phi(\frac{\delta}{2} \sqrt{t})) -(1-\Gamma(t, 1.5, \frac{\delta^2}{8})),
\]
where
\[
\Gamma(t, 1.5, \frac{\delta^2}{8})=\int_0^t \frac{(\delta^2/8)^{3/2} \sqrt{u}}{\Gamma(3/2)} e^{-\delta^2 u/8} dt.
\]
}

From Theorems 1 and 2, we have the following results and the proofs are given in the appendix.

\noindent{\bf Theorem 3.} {\it 

{(i) $E[\sigma_M] = \frac{2}{\delta^2} $ and $Var(\sigma_M) =\frac{12}{\delta^4} $;

(ii) $Cov(\sigma_M, M) = \frac{2}{\delta^3} $  and $\rho (\sigma_M, M) = \frac{1}{\sqrt{3}} $. }
}

The results show that $M$ and $\sigma_M$ are highly correlated.  For this reason, we shall consider to isolate the changed panels mainly based on $M$.

\section{ Approximate null distribution under discrete time model}

We derive the Laplace transform for the joint distribution of $(\sigma_M, M)$.  Let $\tau_+ ^{(0)} =0$. Define
\[
\tau_+  =\tau_+^{(1)}=\inf\{ n>0:  S_n >0 \}
\]
and for $k\geq 2 $
\[
\tau_+^{(k)} =\inf\{n> \tau_+^{(k-1)}: S_n >S_{\tau_+^{(k-1)}} \},
\]
and
\[
K=\sup\{k>0: \tau_+^{(k)} < \infty \},
\]
and $K=0$ if $\tau_+ =\infty $. It can be seen that
\[
P[K=k] =p^k (1-p),
\]
for $k=0,1,2,...$ where $p =P[\tau_+<\infty]$.

Thus, we can write
\[
(\sigma_M , M ) =_d (\tau_+^{(K)}, S_{\tau_+^{(K)}} ).
\]
For given $K=k>0$, $(\tau_+^{(k)}, S_{\tau_+^{(k)}} ) $ is in distribution equivalent to the sum of k i.i.d. copies of $(\tau_+, S_{\tau_+})|\tau_+< \infty$. This leads to the following Laplace transform for $(\sigma_M , M )$ in the normal case.

\noindent{\bf Theorem 4.}  {\it
\[
E[ e^{t\sigma_M +\lambda M}] =\frac{1-G_+ (0,0)}{1-G_+ (s, \lambda)},
\]
where $ G_+ (t, \lambda) =1-E[ e^{t\tau_+ +\lambda S_{\tau_+}}; \tau_+ < \infty] $.
}

\noindent{\it Proof.} 
 By conditioning on the value of $K$, we have
\begin{eqnarray*}
E[ e^{t\sigma_M +\lambda M}] &=& \sum_{k=0}^{\infty} E[ e^{ t \tau_+^{(K)} + \lambda S_{\tau_+^{(K)}}}; K=k] \\
&=& \sum_{k=0}^{\infty} (E[ e^{t \tau_+ + \lambda S_{\tau_+ }}; \tau_+ <\infty])^k P(\tau_+ =\infty ) \\
& =& \frac{1-P(\tau_+ < \infty)}{1- E[ e^{t \tau_+ + \lambda   S_{\tau_+ }}; \tau_+ <\infty]} .
\end{eqnarray*}

 From the above theorem , the exact results for moments of $\sigma_M $ and $M$ can be obtained. For example,
\[
E[\sigma_M ]=\frac{E[\tau_+; \tau_+<\infty]}{P(\tau_+=\infty )}; ~~~ E[M] =\frac{E[S_{\tau_+}; \tau_+<\infty]}{P(\tau_+=\infty )}.
\]

From the Wiener-Hopf factorization (e.g. Siegmund (1985,Theorem 8.41)) , as the random walk $\{S_n\}$ has negative drift, we have
\[
(1-G_+ (t, \lambda))(1-G_- (t, \lambda))=1-e^t e^{\lambda^2/2-\delta \lambda/2},
\]
where $\tau_- =\inf \{n>0: S_n \leq 0 \}$ and $ G_- (t, \lambda)=E[ e^{t \tau_- + \lambda S_{\tau_- }}]$.

By letting $\lambda =0$ and $t\rightarrow 0$, we see
\[
P(\tau_+ =\infty) =\frac{1}{E[\tau_-]}.
\]
 Thus, we have

\noindent{\bf Theorem 5.} 
\[
E[ e^{t\sigma_M +\lambda M}]= \frac{1}{E[\tau_-]} \frac{E[ e^{t \tau_- + \lambda S_{\tau_- }}]}{1-e^{t+\lambda^2/2-\lambda \delta/2}} .
\]

The following corollary shows the second order approximate exponential property for $M$ as $\delta \rightarrow 0$.

\noindent{\bf Corollary 2.}  {\it  As $\delta \rightarrow 0$,
\[
 P( \delta (M+\rho) < x) = 1-e^{- x } (1+o(\delta)) .
\]
 }

\noindent{\it Proof.} 
By taking $t=0$, we have as $\delta \rightarrow 0$,
\begin{eqnarray*}
E[e^{\lambda \delta M}] &=& \frac{1}{E[\tau_-]}\frac{1-E[ e^{ \lambda \delta S_{\tau_- }}]}{1-e^{\lambda^2 \delta^2/2-\lambda \delta^2/2}} \\
&=& \frac{-\delta/2}{E[S_{\tau_-}]} \frac{\lambda \delta E[S_{\tau_-}]+(\delta^2 \lambda^2/2)E[S_{\tau_-}^2]+o(\delta^2)}
{\lambda^2 \delta^2/2 -\lambda \delta^2 /2 +o(\delta^2)} \\
&=& \frac{ 1+ (\lambda \delta /2 ) ( E[S_{\tau_-}^2]/E[S_{\tau_-}]) +o(\delta)}{1-\lambda +o(\delta) } \\
&=& \frac{1}{1-\lambda} e^{\lambda \delta (E[S_{\tau_-}^2]/(2E[S_{\tau_-}]))} +o(\delta^2) .
\end{eqnarray*}
As $\delta \rightarrow 0$, $E[S_{\tau_-}^2]/(2E[S_{\tau_-}]) \rightarrow -\rho $. Thus, we have
\[
E[e^{\lambda \delta (M+\rho )} ] \rightarrow \frac{1}{1-\lambda} .
\]

To study the distribution of $\sigma_M$, we denote by $P^*(.)$ the probability measure with mean $\delta /2$ and $\tau_c =\inf\{ n:  S_n >c \}$ for $c>0$.
We fist note that as $\delta \rightarrow 0$,
\[
P[\sigma_M =0] =P(\tau_+ =\infty) \rightarrow 0.
\]

For a given large value of $M=x$, we use Equation $(3.30)$ of Siegmund (1985) and give the following inverse Gaussian approximation with overshoot correction:
 \[
P[\sigma_M \leq n |M=x] \approx  P^*(\tau_x \leq n )
\]
\[
\approx \Phi (-(x+\rho) n^{-1/2} + (\delta/2) n^{1/2}) + e^{\delta (x+\rho )} \Phi (-(x+\rho) n^{-1/2} - (\delta/2) n^{1/2}).
 \]
In other words, the unconditional distribution of $\sigma_M$ can be treated approximately as a mixture of inverse Gaussian distribution.

\section{Isolation of change panels and estimation of common change point}

Since $M$ and $\sigma_M$ are highly correlated, we shall consider the isolation mainly based on $M$. Conditioning on the common change is detected, we use the corrected exponential distribution for $T_{\tau} (i)$ for $i=1,...N$ for unchanged panels. The BH procedure  (Benjamin and Hochberg (1995)) will be used to control the FDR that is defined as the rate of unchanged panels among all claimed changed panels. Similarly, the FNR is defined as the rate of undiscovered true changed panels among all K true changed panels.

We first calculate the p-values by
\[
p_i = exp(-\delta (T_{\tau}(i) +\rho ))
\]
and $p_{(1)} \leq p_{(2)} \leq .... \leq p_{(N)}$ be the ordered sequence.

For controlled FDR $\alpha $, the number of isolated changed panels will be defined as $\hat{K}$
\[
\hat{K}=\sup \{ i\geq 1: p_{(i)} < \alpha \frac{i}{N} \}.
\]
The well-known theoretical results show that the FDR under this procedure has upper bound $\frac{N-K}{N} \alpha $.
Based on the isolated changed panels, we can estimate the common change point $\nu$ based on the corresponding change point estimations $\hat{\nu}_{(i)}$ for $i=1,..., \hat{K}$.

To show how the proposed procedure performs, we conduct several simulations and leave theoretical investigation for future consideration.

For $\delta =0.5$, $N=100$, $\nu =100$, $ALR_0 =1000$ ($B=74729.5$), and the number of changed panels $K=5,10,20,30$, Table  1 gives the simulation results for FAR (false alarm rate)$P(\tau\leq \nu)$ , FDR, FNR, biases of median estimate $\tilde{\nu}$ and mean estimation $\hat{\nu}$ based on the change-point estimates $\hat{\nu}_i$ from the  $\hat{K}$ isolated changed panels, and mean number $E[\hat{K}]$ of total isolated changed panels $\hat{K}$,  along with the conditional average delay detection time (CADT) $E[\tau-\nu |\tau>\nu]$ based on 5000 simulations. All the values are calculated conditioning on the change is detected $\tau >\nu$.

Table 2 gives the corresponding results for $ARL_0=50000$ ($B=373645.7$) and $\nu=200$.

Figure 1 gives the histograms of simulated FDR, FNR, $\tilde{\nu} -\nu$, and $\hat{K}$ for $ARL=1000$, $\nu=100$, $N=100$, $K=10$, and $\alpha =0.3$ conditioning on $\tau>\nu$.

\begin{figure}
\begin{center}
Figure 1: Histograms of FDR, FNR, $\tilde{\nu}-\nu$, and $\hat{K}$
\includegraphics[width = \textwidth,height=8in]{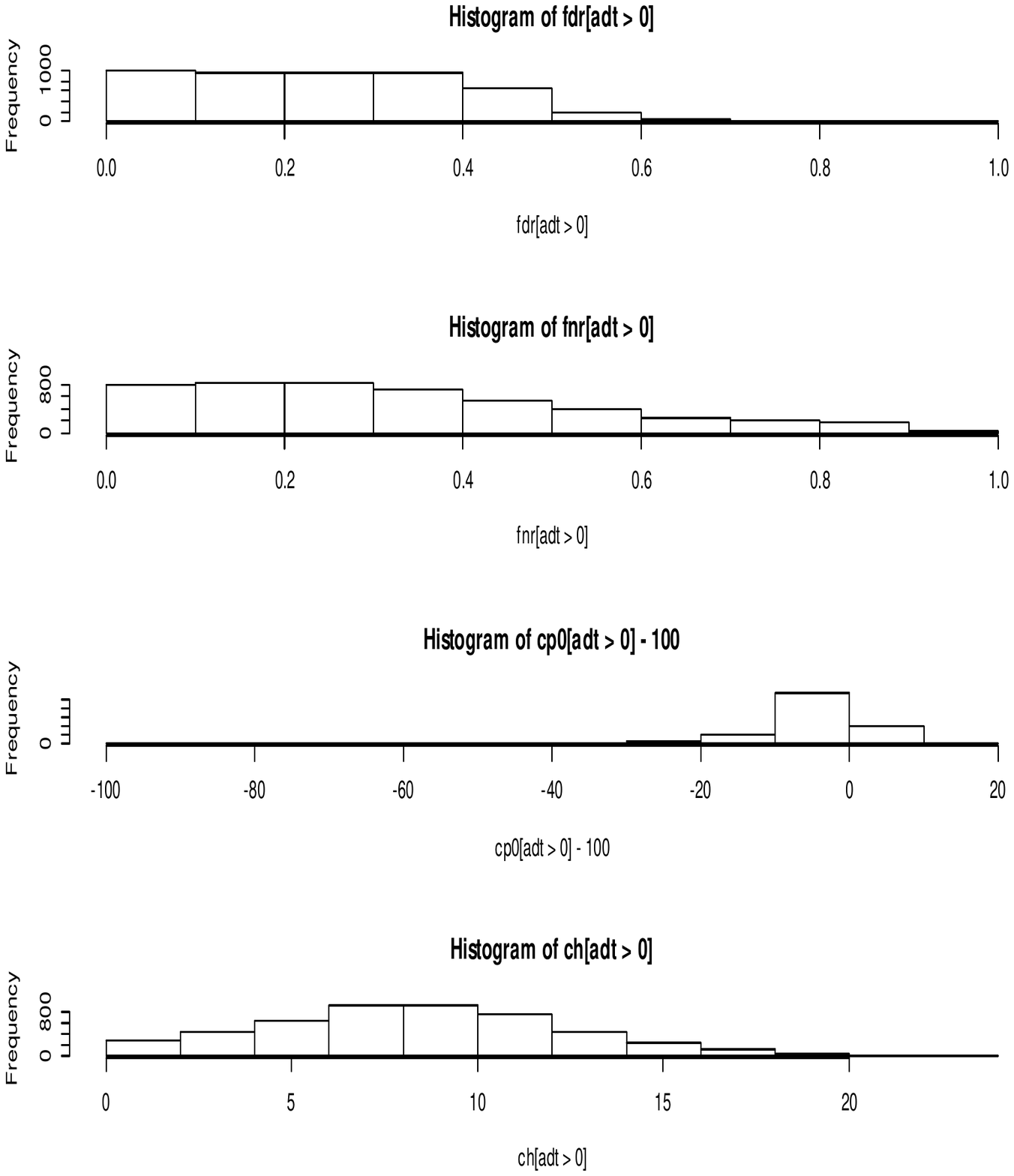}
\end{center}
\end{figure}

\begin{table}
\begin{center}
Table 1. Simulation for $ARL_0=1000$ and $\delta=0.5$ with $N=100$ \\

\begin{tabular}{|c |c | c| c| c| c|}  \hline
K  & $\alpha $ & 0.2 & 0.3 & 0.4 &0.5 \\ \hline
1 & FAR & 0.0372  & 0.042 & 0.0382 & 0.046 \\
  & FDR & 0.216  &0.301  & 0.370 & 0.441 \\
  & FNR & 0.0397  & 0.0401  & 0.0324 & 0.0285 \\
  & $E[\tilde{\nu}-\nu]$ & 0.0 &  1.0 & 3.0 & 5.0 \\
  &$E[\hat{\nu}- \nu]$  & 1.368  & 4.01 & 6.52 & 9.88 \\
  & $E[\hat{K}]$ &  1.481      & 1.817  & 2.299 & 2.935 \\
  & CADT &    60.03   & 59.70 & 59.51 & 60.27 \\ \hline
5 & FAR & 0.0406 & 0.0446 & 0.0422 & 0.0488 \\
  & FDR & 0.188 & 0.270 & 0.349 & 0.431 \\
  & FNR & 0.417& 0. 346 & 0.282 & 0.243 \\
  & $E[\tilde{\nu}-\nu]$ & -3 & -2 & -1& 0 \\
  & $E[\hat{\nu}- \nu]$ & -6.6 & -4.6 & -2.9 & -1 \\
  & $E[\hat{K}]$ & 3.83 & 4.99 & 6.38 & 8.04 \\
  & CADT & 32.89  & 32.98 &33.21 & 33.01 \\ \hline
10&FAR &0.0424 & 0.0398 & 0.0434 & 0.0428 \\
  & FDR & 0.172 & 0.256 & 0.333 & 0.422 \\
  & FNR & 0.469 & 0.375 & 0.302 & 0.250 \\
  & $E[\tilde{\nu}-\nu]$& -3 & -2 & -1 & 0 \\
  & $E[\hat{\nu}- \nu]$ &-6.46 & -5.0 & -3.20 & -1.6 \\
  & $E[\hat{K}]$ & 6.65 & 8.88 & 11.31 & 14.36 \\
  & CADT & 26.31 & 26.10 & 26.40 & 26.50 \\ \hline
20& FAR & 0.0432 & 0.0404 & 0.0418 & 0.0388 \\
  & FDR & 0.155 & 0.227 & 0.307 & 0.383 \\
  & FNR & 0.482 & 0.374 & 0.283 & 0.216 \\
  & $E[\tilde{\nu}\nu]$ &-3 & -2 & -1 & 0 \\
  & $E[\hat{\nu}- \nu]$ &-6.4 & -4.5 & -2.9 & -1.2 \\
  & $E[\hat{K}]$& 12.47 & 16.7 & 21.5 & 26.5 \\
  &CADT& 21.52 & 21.30 & 21.42 & 21.29 \\ \hline
30 & FAR & 0.047 & 0.0438 & 0.0428 & 0.046 \\
  & FDR & 0.140 & 0.205 & 0.273 & 0.342 \\
  & FNR& 0.480 & 0.348 & 0.265 & 0.192 \\
  & $E[\tilde{\nu}-\nu]$  &-3 & -2 & -1 & 0 \\
  & $E[\hat{\nu}- \nu]$ & -6.19 & -4.1 & -2.8 & -1.4 \\
  & $E[\hat{K}]$ & 18.34 & 25.0 & 30.9 & 37.57 \\
  &CADT & 18.29 & 18.57 & 18.44 & 18.65 \\ \hline

\end{tabular}
\end{center}
\end{table}

\begin{table}
\begin{center}
Table 2. Simulation for $ARL_0=5000$ and $\delta=0.5$ with $N=100$ \\

\begin{tabular}{|c |c | c| c| c| c|}  \hline
K  & $\alpha $ & 0.2 & 0.3 & 0.4 &0.5 \\ \hline
1 & FAR & 0.0232  & 0.0284 & 0.0272& 0.0232 \\
  & FDR & 0.184  &0.277  & 0.349 & 0.427 \\
  & FNR & 0.0082  & 0.0082  & 0.0047 & 0.0055 \\
  & $E[\tilde{\nu}-\nu]$ & 1.0 &  3.0 & 6.0 & 10.0 \\
  &$E[\hat{\nu}- \nu]$  & 4.58  & 8.37 & 11.68 & 17.20 \\
  & $E[\hat{K}]$ &  1.468      & 1.853  & 2.299 & 2.960 \\
  & CADT &    75.07  & 75.21 & 74.47 & 74.94 \\ \hline
 5 & FAR & 0.0244 & 0.0244 & 0.0278 & 0.0282 \\
  & FDR & 0.185 & 0.268 & 0.351 & 0.434 \\
  & FNR & 0.265& 0.214 & 0.176 & 0.140 \\
  &  $E[\tilde{\nu}-\nu]$ & -1 & 0 & 1& 2 \\
  & $E[\hat{\nu}- \nu]$ & -2.7 & -1.0 & 0.9 & 3.9 \\
  & $E[\hat{K}]$ & 4.79 & 5.91 & 7.21 & 8.99 \\
  & CADT & 42.84  & 42.96 &42.95 & 42.90 \\ \hline
10&FAR &0.0292 & 0.0224 & 0.0212 & 0.0242 \\
  & FDR & 0.172 & 0.256 & 0.338 & 0.427 \\
  & FNR & 0.296 & 0.224 & 0.176 & 0.139 \\
  & $E[\tilde{\nu}-\nu]$ & -1 & 0 & 0 & 1.5 \\
  & $E[\hat{\nu}- \nu]$ &-2.9 & -1.6 & 0.3 & 2.6 \\
  & $E[\hat{K}]$ & 8.78 & 10.94 & 13.3 & 16.35 \\
  & CADT & 35.47 & 35.37 & 35.40 & 35.46 \\ \hline
20& FAR & 0.0234 & 0.0280 & 0.0234 & 0.0244 \\
  & FDR & 0.158 & 0.234 & 0.314 & 0.388 \\
  & FNR & 0.280 & 0.211 & 0.158 & 0.115 \\
  & $E[\tilde{\nu}-\nu]$ &-1 & -0.5 & 0 & 1 \\
  & $E[\hat{\nu}- \nu]$ &-2.4 & -1.3 & -0.1 & 1.6 \\
  & $E[\hat{K}]$ & 17.34 & 21.04 & 25.2 & 29.9 \\
  &CADT& 29.55 & 29.41 & 29.41 & 29.54 \\ \hline
30 & FAR & 0.0218 & 0.0242 & 0.0248 & 0.0252 \\
  & FDR & 0.139 & 0.207 & 0.274 & 0.344 \\
  & FNR& 0.271 & 0.193 & 0.138 & 0.102 \\
  & $E[\tilde{\nu}-\nu]$ &-1 & 0 & 0 & 1 \\
  & $E[\hat{\nu}- \nu]$ & -2.1 & -1.2 & -0.1 & 0.98 \\
  & $E[\hat{K}]$ & 25.63 & 30.90 & 36.17 & 41.80 \\
  &CADT & 26.38 & 26.52 & 26.59 & 26.48 \\ \hline

\end{tabular}
\end{center}
\end{table}

By looking at the simulation results of Tables 1 and 2, we have several important findings.

(i) The FDRs are not significantly different between $ARL_0 =1000 $ and $ARL_0 =5000$ and decreases when K increases;

(ii) The FNRs are not significantly different when K changes for fixed $ARL_0$ and decreases when $ARL_0$ increases;

(iii) The simulated FDRs are very close to the theoretical upper bound $ \alpha \frac{N-K}{N} $;

(iv) The FNRs decreases as $\alpha $ increases and the two are roughly balanced around $\alpha =0.3$. So $\alpha =0.2 $ or 0.3 are recommended.

(v) The median estimate for the common change point is preferred as its bias are smaller than the mean estimate;

(vi) The number of isolated panels increases as $\alpha $ increases and roughly equals to the true K at $\alpha =0.3$.

However as the post-change mean $\mu$ is rarely known and we typically select $\delta$ as the minimum magnitude to detect. So we also run a simulation study for $\mu =1.0, 1.5$ and 2.0. Table 3 gives the results for both $ARL_0 =1000$ and $5000$ with $\alpha =0.2$ and 0.3. Additional findings are:

(vii) The FDRs are roughly the same when $\mu $ changes, while the FNRs are reduced more significantly for $\mu > \delta $ as $\mu$ increases.
     From this point of view, $\alpha =0.2$ is preferred if stronger signals are expected.

(viii) However, as $\mu$ increases, the bias of the common change point becomes more negative, similar to Table 2.1 in Wu (2005, pg 40).

\begin{table}
\begin{center}
Table 3. Simulation for unknown $\mu $  \\

\begin{tabular}{|c |c | c| c| c| c|}  \hline
  & $\mu $ & 0.5 & 1.0 & 1.5 & 2.0 \\ \hline
$ARL_0 =1000$ & FAR & 0.0398 & 0.0392 & 0.0420 & 0.0432 \\
$\nu=100$ & FDR & 0.256 & 0.254 & 0.261 & 0.259 \\
$K=10$  & FNR & 0.375& 0. 152 & 0.074 & 0.040 \\
 $\alpha =0.3$ &$E[\tilde{\nu}-\nu]$ & -2 & -4 & -5& -5.5 \\
  & $E[\hat{\nu}- \nu]$  & -5.0 & -6.42 & -6.74 & -6.88 \\
  & $E[\hat{K}]$& 8.88 & 11.86 & 13.01 & 13.40 \\
  & CADT & 26.5  & 12.14 & 8.06 & 6.08 \\ \hline
$ARL_0 =1000$&FAR &0.0424 & 0.0442 & 0.0414 & 0.0368 \\
$\nu=100$   & FDR & 0.172 & 0.173 & 0.172 & 0.172 \\
$K=10$  & FNR & 0.469 & 0.226 & 0.128 & 0.080 \\
$\alpha =0.2$  & $E[\tilde{\nu}-\nu]$ & -3 & -4.5 & -5 & -4.5 \\
  & $E[\hat{\nu}- \nu]$  &-6.46 & -7.28 & -6.89 & -6.71 \\
  & $E[\hat{K}]$ & 6.65 & 9.63 & 10.77 & 11.38 \\
  & CADT & 26.31 & 12.21 & 8.05 & 6.07 \\ \hline
$ARL_0 =5000$& FAR & 0.0224 & 0.0278 & 0.0254 & 0.0230 \\
$\nu=200$  & FDR & 0.256 & 0.258 & 0.265 & 0.268 \\
$K=10$  & FNR & 0.224 & 0.048 & 0.0167 & 0.007 \\
$\alpha =0.3$  & $E[\tilde{\nu}-\nu]$  &0 & -3& -4.5 & -5 \\
  & $E[\hat{\nu}- \nu]$ &-1.6 & -4.4 & -5.5 & -5.9 \\
  & $E[\hat{K}]$ & 10.94 & 13.3 & 13.9 & 14.01 \\
  &CADT& 35.37 & 15.95 & 10.42 & 7.82 \\ \hline
$ARL_0 =5000$ & FAR & 0.0292 & 0.0258 & 0.0262 & 0.0294 \\
$\nu=200$  & FDR & 0.172 & 0.172 & 0.176 & 0.178 \\
 $K=10$  & FNR& 0.296 & 0.084 & 0.035 & 0.015 \\
$\alpha =0.2$  &  $E[\tilde{\nu}-\nu]$  &-1 & -3 & -4 & -4\\
  & $E[\hat{\nu}- \nu]$ & -2.9 & -4.8 & -5.2 & -5.4 \\
  & $E[\hat{K}]$ & 8.78 & 11.33 & 11.97 & 12.24 \\
  &CADT & 35.47 & 15.93 & 10.41 & 7.82 \\ \hline

\end{tabular}
\end{center}
\end{table}

\section{Conclusion}

In this paper, we proposed a BH procedure to control the FDR after a common change is detected in multi-panel data stream. The method only uses partial information available from each individual CUSUM process and is shown performing quire well.
To reduce the FDR for isolating changed panels and estimating the common change point, supplementary runs are necessary on isolated changed panels. A simple method is to run one-sided truncated sequential tests by just finding the true changed panels as discussed in Wu (2018).
Further discussions on sequential multiple tests on controlling FDR can also be used in the supplementary runs; see Bartroff (2017), De and Baron (2015), and Song and Fellouris (2019). As discussed in Wu (2019), we may also use the adaptive combined SR-CUSUM procedure which can eliminate large biases of the common change point estimation when the post change means are unknown. The results will be presented in future communications.

\section {Appendix}

\subsection{Proof of Theorem 1}

\begin{eqnarray*}
& & P[\sigma_M<t, M>x] \\
& = & P[\sup_{s\geq t} W_s < \sup_{0\leq s < t} W_s <x ] \\
&= & P[W_t +\sup_{0\leq s <\infty} W_s' < \sup_{0\leq s<t} W_s < x ] \\
&=& P[\sup_{0\leq s<t} W_s < x , W_t<x-M'] -P[\sup_{0\leq s<t} W_s < W_t+M', W_t<x-M'] \\
&=& P[\sup_{0\leq s<t} W_s < x , W_t<x-M']-P^*[\sup_{0\leq s<t} W_s <M', W_t >M'-x] \\
&=&P^*[\sup_{0\leq s<t} W_s >M', W_t >M'-x]-P[\sup_{0\leq s<t} W_s > x , W_t<x-M'] \\
&=& P^*[\sup_{0\leq s<t} W_s >M']-P^*[\sup_{0\leq s<t} W_s >M', W_t <M'-x] \\
& & -P[\sup_{0\leq s<t} W_s > x , W_t<x-M'], \\
\end{eqnarray*}
where in the second equation from last, we use the fact
\[
P[W_t >M' -x] = P^*[W_t <x-M' ] .
\]

The last two terms  are evaluated by using Equation (3.14) of Siegmund(1985):
 \[
 P^*[\sup_{0\leq s<t} W_s >M', W_t <M'-x] =E[e^{\delta M'} \Phi (-\frac{x+M'}{\sqrt{t}} -\frac{\delta}{2}\sqrt{t} )];
 \]
 \[
 P[\sup_{0\leq s< t} W_s > x , W_t<x-M'] =E[e^{\delta x} \Phi (-\frac{x+M'}{\sqrt{t}} +\frac{\delta}{2}\sqrt{t} )],
 \]
 where $\Phi(x) $ and $\phi(x)$ are standard normal cdf and pdf.

\subsection{Proof of Theorem 2}

\begin{eqnarray*}
& &E[ \Phi (-\frac{x+M'}{\sqrt{t}} +\frac{\delta}{2}\sqrt{t} ) +\frac{2}{\delta \sqrt{t}} \phi (-\frac{x+M'}{\sqrt{t}}+\frac{\delta}{2}\sqrt{t} )]
\\
&& =\int_0^{\infty} \delta e^{-\delta y} \Phi (-\frac{x+y}{\sqrt{t}} +\frac{\delta}{2}\sqrt{t} ) dy
+\frac{2}{\sigma \sqrt{t}} \int_0^{\infty} \delta e^{-\delta y} \phi (-\frac{x+y}{\sqrt{t}}+\frac{\delta}{2}\sqrt{t} )dy \\
&&= \Phi(-\frac{x}{\sqrt{t}}+\frac{\delta}{2}\sqrt{t} )
+\int_0^{\infty} (-\frac{1}{\sqrt{t}}) e^{-\delta y} \phi (\frac{x+y}{\sqrt{t}}
-\frac{\delta}{2}\sqrt{t} ) dy \\
& & +\frac{2}{\sqrt{t}} \int_0^{\infty} e^{-\delta y} \phi (\frac{x+y}{\sqrt{t}}-\frac{\delta}{2}\sqrt{t} )dy \\
&&=\Phi(-\frac{x}{\sqrt{t}}+\frac{\delta}{2}\sqrt{t} ) +\frac{1}{\sqrt{t}} \int_0^{\infty} e^{-\delta y} \phi (\frac{x+y}{\sqrt{t}}-\frac{\delta}{2}\sqrt{t} )dy \\
&&= \Phi(-\frac{x}{\sqrt{t}}+\frac{\delta}{2}\sqrt{t} ) + \frac{1}{\sqrt{t}} \int_0^{\infty} e^{-\delta y} \frac{1}{\sqrt{2\pi}} e^{-(y+x-(\delta/2)t)^2/(2t)} dy \\
&&= \Phi(-\frac{x}{\sqrt{t}}+\frac{\delta}{2}\sqrt{t} ) +\frac{1}{\sqrt{t}} \int_0^{\infty}\frac{1}{\sqrt{2\pi}} e^{-(y+(x+\delta t/2))^2/(2t)} dy e^{\delta x} \\
&&= \Phi(-\frac{x}{\sqrt{t}}+\frac{\delta}{2}\sqrt{t} ) +e^{\delta x} (1- \Phi (\frac{x}{\sqrt{t}}+\frac{\delta}{2}\sqrt{t} )) \\
&&= \Phi(-\frac{x}{\sqrt{t}}+\frac{\delta}{2}\sqrt{t} ) +e^{\delta x}\Phi (-\frac{x}{\sqrt{t}}-\frac{\delta}{2}\sqrt{t} ).\\
\end{eqnarray*}

\subsection{Proof of Theorem 3}

First, we note $E[M] =\frac{1}{\delta}$ and $Var(M) =\frac{1}{\delta^2}$. Second, by using the property of inverse Gaussian distribution,
\[
E[\sigma_M] =E[E[\sigma_M|M]] =E[\frac{M}{\delta/2}] = \frac{2}{\delta^2} ,
\]
\begin{eqnarray*}
Var(\sigma_M) &=& E[Var(\sigma_M|M)] +var(E[\sigma_M|M]) \\
&=& E[\frac{M}{(\delta/2)^3}] +Var(\frac{M}{\delta/2}) \\
&=& \frac{1}{\delta(\delta/2)^3} +\frac{1}{\delta^2 (\delta/2)^2} = \frac{12}{\delta^4} \\
\end{eqnarray*}
Also,
\[
E[\sigma_M M] =E[ME[\sigma_M|M]] =E[M^2/(\delta/2)] =\frac{4}{\delta^3} .
\]
(ii) is proved by combining the above results.\\

\noindent{\bf Acknowledgement.} 

This research is partially supported by a RSCA grant from California State University at Stanislaus. \\

\noindent{\bf References}
{\small
\begin{description}
\item[{\rm Bartroff, J., 2018.}] Multiple hypothesis tests controlling generalized error rates
for sequential data. {\em Statistica Sinica} \textbf{28}, 363--398.
\item[{\rm Benjamini,Y., Hochberg,Y.,1995.}] Controlling
the false discovery rate: A practical and powerful approach to
multiple testing. {\em J. Roy. Statist. Soc. (B)}
\textbf{57}, 289 -- 300.
\item[{\rm Chan, H. P.,2017.}] Optimal sequential detection in multi-stream data. {\em Annals of Statistics} \textbf{45}(6),2636--2763.
\item[{\rm De, S. K., Baron, M.,2015.}] Sequential tests controlling generalized familywise
error rates. {\em Statistical Methodology} \textbf{23},88 -- 102.
\item[{\rm Mei, Y.,2010.}] Efficient scalable schemes for monitoring a large number of data streams. {\em Biometrika} \textbf{97},419 --433.
\item[{\rm Pollak, M., 1987.}] Average run lengths of an optimal method for detecting a change in distribution. {\em Annals of Statistics} \textbf{15},749 --779.
\item[{\rm Siegmund, D., 1985.}] {\em Sequential Analysis: Tests and Confidence Intervals}. New York: Springer.
\item[{\rm Song, Y., Fellouris, G., 2019.}] Sequential multiple testing with generalized error control: An asymptotic optimality theory.  {\em Ann. Statist} \textbf{47(3)}, 1776 -– 1803.
\item[{\rm Tartakovsky, A.G., Veeravalli, V.V., 2008.}] Asymptotically optimal quickest detection change detection in distributed sensor. {\em Sequential Analysis } \textbf{27}, 441 --475.
\item[{\rm Wu, Y., 2005.}] {\em Inference for Change-point and Post-change Means After a CUSUM Test.}
{\em Lecture Notes in Statistics } 180, Springer, New York .
\item[{\rm Wu, Y., 2018.}] Supplementary score test for sparse signals in large-scale truncated sequential tests. {\em Journal of Statistical Theory and Practice} \textbf{12}(4),744--756.
\item[{\rm Wu, Y., 2019.}] A combined SR-CUSUM procedure for detecting common changes in panel data. {\em Communication in Statistics: Theory and Methods} \textbf{48}(17): 4302--4319.
\item[{\rm Xie, Y., Siegmund, D.,2013.}] Sequential multi-sensor change-point detection. {\em Annals of Statistics} \textbf{41},670 --692.

\end{description}
}

\end{document}